\documentclass[letterpaper, 10 pt, conference]{ieeeconf}  

\IEEEoverridecommandlockouts                              

\usepackage{graphics} 
\usepackage{epsfig} 
\usepackage{mathptmx} 
\usepackage{times} 
\usepackage{amsmath} 
\usepackage{amssymb}  
\usepackage{amsbsy,bm,fixmath}

\newtheorem{theorem}{Theorem}
\newtheorem{definition}[theorem]{Definition}
\newtheorem{lemma}[theorem]{Lemma}
\newtheorem{proposition}[theorem]{Proposition}

\newtheorem{remark}{Remark}

\def\mbA{\mathbold{A}}
\def\mbB{\mathbold{B}}
\def\mbD{\mathbold{D}}
\def\mbP{\mathbold{P}}
\def\mbx{\mathbold{x}}
\def\mbQ{\mathbold{Q}}
\def\mbS{\mathbold{S}}
\def\mbK{\mathbold{K}}

\def\mbr{\mathbold{r}}

\title{\LARGE \bf Linear Quadratic Mean Field Games -- Part I: The Asymptotic\\ Solvability Problem}

\author{Minyi Huang\qquad Mengjie Zhou
\thanks{This work was  supported by Natural Sciences and
Engineering Research Council (NSERC) of Canada. In Proc. 23rd MTNS,
 2018, pp. 489-495.}
\thanks{M. Huang is with the School of
Mathematics and Statistics, Carleton University, Ottawa, K1S 5B6 ON,
Canada (mhuang@math.carleton.ca).  }%
\thanks{M. Zhou is with the School of Mathematics and Statistics, Carleton University, Ottawa, K1S 5B6 ON,
Canada (mengjiezhou@cmail.carleton.ca).  }
}

\begin{document}

\maketitle
\thispagestyle{empty}
\pagestyle{empty}

\begin{abstract}
This paper
investigates the so-called asymptotic solvability problem in linear quadratic (LQ) mean field games.
The model has asymptotic solvability if for all sufficiently large population sizes, the corresponding game has a set of feedback Nash strategies subject to a mild regularity requirement.
We provide a necessary and sufficient condition and show that in this case the solution converges to a mean field limit. This is accomplished by developing a re-scaling method to derive a low dimensional ordinary differential equation (ODE) system, where a non-symmetric Riccati ODE has a central role.
\end{abstract}




\section{INTRODUCTION}

Mean field game theory has undergone a phenomenal growth. It provides a powerful methodology
for handling complexity in noncooperative mean field decision problems.
The readers are referred to \cite{BFY13,CHM17} for an overview of the theory and applications.
 The analysis in the LQ setting has attracted
substantial interest due to its appealing analytical structure \cite{HCM07,LZ08,MB15}. Specifically, the strategy of an individual player can be determined in a feedback form using its own state.

Two important methodologies called the top-down approach and
bottom-up approaches \cite{CHM17}, respectively,
have been widely used in the analysis of mean field games.
By the top-down approach \cite{HCM07,HMC06}, one determines the best response of a representative agent to a mean field of an infinite population, and next all the agents's best responses should regenerate that mean field.
This procedure formalizes a fixed point problem which can be solved and further used to design decentralized strategies. The bottom-up
 approach (also called the direct approach \cite{HZ18})
   starts by formally solving an $N$-player game to obtain a large coupled solution equation system. The next step is to derive a simple liming equation system by taking $N\to \infty$ \cite{LL07}; also see \cite{L16} for a probabilistic framework.

This paper considers the LQ mean field game and addresses the so-called asymptotic solvability.
We start with an entirely conventional solution of the game by dynamic programming and
derive a set of coupled Riccati ODEs. This method may be viewed as an instance of the bottom-up approach.
 Our objective is to find a necessary and sufficient condition
for the sequence of games to be appropriately solvable. It turns out that such a condition
is completely determined by  a single low dimensional non-symmetric Riccati ODE.
The derivation of this condition involves a novel re-scaling method for large-scale coupled equations with two-scales of interactions. We further determine the  mean field limit of the individual strategies.
Our approach has  connection with an early model of mean field social optimization, which studies a high dimensional algebraic Riccati equation and uses symmetry for dimension reduction \cite[Sec. 6.3]{H03}.
 The methodology of extracting a low dimensional structure, here as a non-symmetric Riccati ODE, to capture essential information of the large scale decision model shares similarity to identifying low dimensional dynamics of coupled oscillators in the physics literature \cite{MBS09,OA08,PM14}.
Other related works include \cite{H15,P14,P15}.
An optimal control problem for a set of agents with mean field coupling
is solved in \cite{H15} by a large-scale Riccati ODE, where a mean field limit is derived for the Riccati equation using the scalar state and symmetry.
An LQ Nash game of infinite horizon is analyzed in \cite{P14} where the number of players increases to infinity. The method is to postulate the strategies of all players and  examine the control problem of a fixed player subject to the mean field dynamics. Then a family of low dimensional control problems and the  parameterized algebraic Riccati equations  can be solved by an implicit function theorem. Sufficient conditions are obtained for solvability when the population size is large. The solvability of LQ games with increasing population sizes in the set-up of \cite{LL07} is studied in \cite{P15} analyzing $2N$-coupled steady-state Hamilton-Jacobi-Bellman (HJB) and Fokker-Planck-Kolmogorov (FPK) equations, where each player's control is restricted to be local state feedback from the beginning. Some algebraic conditions are obtained. However, it requires some restrictions on the model parameters, including symmetric state coefficients in dynamics.

Note that the top-down approach can also be  used to solve the LQ mean field game \cite{CHM17}.
In part II \cite{HZ18} of this paper, we  investigate the relation between the top-down approach and
the bottom-up approach as developed in this paper. A surprising finding is  that they are not
always equivalent.

The organization of the paper is as follows.
Section \ref{sec:mod} describes the LQ Nash game together with its solution via dynamic programming and Riccati equations.
Section \ref{sec:asy} presents the necessary and sufficient condition for asymptotic  solvability, for which we give the proof in Section \ref{sec:proof}. Section \ref{sec:dec} presents further mean field limits related to the dynamic programming equation and derives decentralized strategies. An illustrative example is provided in Section \ref{sec:exa}.
Section \ref{sec:con} concludes the paper.

\section{The LQ Nash Game} \label{sec:mod}

Consider a population of $N$ players (or agents) denoted by ${\cal A}_i$,  $1\le i\le N$.
The state process $X_i(t)$ of ${\cal A}_i$ satisfies the following
stochastic differential equation (SDE)
\begin{align}\label{stateXi}
dX_i(t)=\big(AX_i(t)+Bu_i(t)+GX^{(N)}(t)\big)dt+DdW_i(t),
\end{align}
where the state  $X_i\in \mathbb{R}^n$, control
$u_i\in\mathbb{R}^{n_1}$, and
 $ X^{(N)}=\frac{1}{N}\sum_{k=1}^N X_k$. The initial states $\{X_i(0), 1\le i\le N\}$ are independent with $EX_i(0)=x_i(0)$ and finite second moment. The $N$ standard  $n_2$-dimensional Brownian motions $\{W_i, 1\le i\le N\}$ are independent and also independent of the  initial states.
For symmetric matrix $S\ge 0$, we may write $x^TSx=|x|_S^2$.     
The cost of player ${\cal A}_i$  is given by
\begin{align}\label{costJi}
J_i =\  & E\int_0^T \Big(  |X_i(t)-\Gamma X^{(N)}(t)-\eta|_Q^2+u_i^T(t) R u_i(t)\Big)dt \nonumber\\
&+E|X_i(T)-\Gamma_f X^{(N)}(T)-\eta_f|_{Q_f}^2.
\end{align}
The constant matrices (or vectors)
$A$, $B$, $G$, $D$, $\Gamma$, $Q$, $R$, $\Gamma_f$, $Q_f$, $\eta$, $\eta_f$ above
have compatible dimensions, and $Q\ge 0$, $R>0$, $Q_f\ge 0$.
For notational simplicity, we only consider constant parameters
for the model. Our analysis and results can be easily extended to
the case of time-dependent parameters.


Define
\begin{align}
&X(t)= \begin{bmatrix} X_1(t) \\ \vdots \\ X_N(t)
\end{bmatrix}\in \mathbb{R}^{Nn},\nonumber
\quad W(t)=\begin{bmatrix}W_1(t) \\ \vdots \\ W_N(t)
\end{bmatrix}\in \mathbb{R}^{Nn_2},\nonumber   \\
&\widehat \mbA = \mbox{diag}[A, \cdots, A]+{\bf 1}_{n\times n}\otimes
\frac{G}{N}\in\mathbb{R}^{Nn\times Nn},\nonumber \quad \\
&\widehat \mbD =\mbox{diag}[D, \cdots, D]\in \mathbb{R}^{Nn\times Nn_2},\nonumber \\
&\mbB_k = e_k^N\otimes B\in\mathbb{R}^{Nn\times n_1}, \quad 1\le k\le N. \nonumber \quad
\end{align}
We denote by ${\bf 1}_{k\times l}$  a $k\times l$ matrix with all entries equal to 1, by $\otimes$ the Kronecker product, and by the column  vectors  $\{e_1^k, \ldots, e_k^k\}$  the canonical basis of $\mathbb {R}^k$.  We may use a subscript $n$ to indicate the identity matrix $I_n$ to be $n\times n$.

Now we write \eqref{stateXi} in the form
\begin{align}\label{bigx}
dX(t)=\Big(\widehat \mbA X(t)+\sum_{k=1}^N \mbB_k u_k(t)\Big)dt+\widehat \mbD dW(t).
\end{align}
Under closed-loop state information, we denote the value function of ${\cal A}_i$  by $V_i(t,\mbx)$, $1\le i\le N$, which corresponds to the initial condition $X(t)= \mbx=(x_1^T, \ldots, x_N^T)^T$.
The set of value functions is determined by the system of  HJB equations
\begin{align}\label{DE}
&0= \frac{\partial V_i}{\partial t}+\min_{u_i\in {\mathbb R}^{n_1}}\Bigg(\frac{\partial^T V_i}{\partial \mbx }\big(\widehat {\mbA}\mbx+\sum_{k=1}^N {\mbB}_k u_k\big)
+u_i^T R u_i \nonumber \\
   &\qquad+|x_i-\Gamma x^{(N)}-\eta|_Q^2  +\frac{1}{2}\mbox{Tr}\Big
   ({\widehat{\mbD}^T (V_i)_{\mbx\mbx} \widehat{\mbD}}\Big)\Bigg), \\
&V_i(T,\mbx)=|x_i-\Gamma_f x^{(N)}-\eta_f|_{Q_f}^2, \quad  1\le i\le N, \nonumber
\end{align}
where $x^{(N)}=(1/N)\sum_{k=1}^N x_k$ and the minimizer is
 $$u_i=-\frac{1}{2} R^{-1} {\mbB}_i^T \frac{\partial V_i}
 {\partial \mbx}.$$  Next we substitute $u_i$ into \eqref{DE}:
\begin{align}\label{DE2}
0 = &\frac{\partial V_i}{\partial t}+\frac{\partial^T V_i}{\partial \mbx}\Big({\widehat \mbA}\mbx
-\sum_{k=1}^N \frac{1}{2} {\mbB}_k R^{-1} {\mbB}_k^T \frac{\partial V_k}{\partial \mbx}\Big) \nonumber \\
&+|x_i-\Gamma x^{(N)}-\eta|_Q^2  \nonumber \\
   &+\frac{1}{4}\frac{\partial^T V_i}{\partial \mbx}{\mbB}_i R^{-1} \mbB_i^T\frac{\partial V_i}{\partial \mbx}+\frac{1}{2}\mbox{Tr}\Big({{\widehat \mbD}^T (V_i)_{\mbx\mbx} {\widehat \mbD}}\Big).
\end{align}

Denote
\begin{align*}
&\mbK_i =[0,\cdots,0, {I}_n,0,\cdots,0]-\frac{1}{N}
[\Gamma,\Gamma,\cdots,\Gamma],\\
& \mbK_{if} =[0,\cdots,0, {I}_n,0,\cdots,0]-\frac{1}{N}
[\Gamma_f,\Gamma_f,\cdots,\Gamma_f],\\
& \mbQ_i=\mbK_i^T Q \mbK_i, \quad  \mbQ_{if}=\mbK_{if}^T Q_f \mbK_{if},
\end{align*}
where $I_n$ is the $i$th submatrix.
We write
\begin{align} 
|x_i-\Gamma x^{(N)}-\eta|_Q^2  = &\mbx^T {\mbQ}_i \mbx - 2\mbx^T
\mbK_i^TQ\eta + \eta^T Q\eta, \nonumber
    \end{align}
and we can write $|x_i-\Gamma_f x^{(N)}-\eta|_{Q_f}^2 $ in a similar form.

Suppose $V_i(t,\mbx)$ has the following form
\begin{align}\label{Vform}
&V_i(t,\mbx)=\mbx^T {\mbP}_i(t) \mbx+2{\mbS}_i^T(t) \mbx+\mbr_i(t),
\end{align}
where $\mbP_i$ is symmetric.
Then
\begin{align}  \label{1st}
\frac{\partial V_i}{\partial \mbx}=2{\mbP}_i(t)\mbx+2{\mbS}_i(t),\quad
\frac{\partial^2 V_i}{\partial \mbx^2}=2\mbP_i(t).
\end{align}

We substitute \eqref{Vform} and \eqref{1st}  into \eqref{DE2} and derive the equation systems:
\begin{align}\label{DE3_P}
\begin{cases}
\dot{\mbP}_i(t) =  - \Big({\mbP}_i(t)\widehat{\mbA}+\widehat{\mbA}^T
\mathbb{\mbP}_i(t)\Big)+\\
                          \qquad\qquad     \Big({\mbP}_i(t)\sum_{k=1}^N
                             {\mbB}_k R^{-1} {{\mbB}_k}^T {\mbP}_k(t)\\
                      \qquad\qquad         +\sum_{k=1}^N {\mbP}_k(t){\mbB}_k R^{-1} {\mbB}^T_k {\mbP}_i(t)\Big) \\
                         \qquad\qquad      - {\mbP}_i(t){\mbB}_i R^{-1} {\mbB}_i^T
                             {\mbP}_i(t)- {\mbQ}_i   , \\
 {\mbP}_i(T) = {\mbQ}_{if},
 \end{cases}
\end{align}

\begin{align}\label{DE3_S}
\begin{cases}
\dot{{\mbS}}_i(t) = - {\widehat\mbA}^T {\mbS}_i(t)  - {\mbP}_i(t){\mbB}_i R^{-1} {{\mbB}_i}^T {\mbS}_i(t)\\
                            \qquad\qquad    +{\mbP}_i(t) \sum_{k=1}^N {\mbB}_k R^{-1} {\mbB}_k^T {\mbS}_k(t)\\
                             \qquad\qquad   + \sum_{k=1}^N {\mbP}_k(t){\mbB}_k R^{-1}{\mbB}^T_k {\mbS}_i(t)  \\
       \qquad\qquad   +\mbK_i^T Q\eta , \\
{\mbS}_i(T)= -\mbK_{if}^T Q_f\eta_f,
\end{cases}
\end{align}

\begin{align}\label{DE3_gamma}
\begin{cases}
\dot{\mbr}_i(t) =  2{{\mbS}_i}^T(t)\sum_{k=1}^N {\mbB}_k R^{-1}
{\mbB}_k^T {\mbS}_k(t) \\
                   \qquad\qquad  - {\mbS}_i^T(t){\mbB}_i R^{-1} {\mbB}_i^T{\mbS}_i(t)\\
                    \qquad\qquad         - \eta^T Q \eta   -\mbox{Tr}\big(\widehat{\mbD}^T
                          {\mbP}_i(t)\widehat{\mbD}\big), \\
 \mbr_i(T)=\eta_f ^T Q_f \eta_f.
\end{cases}
\end{align}

\begin{remark}\label{remark:P}
If \eqref{DE3_P} has a solution $(\mbP_1, \cdots, \mbP_N)$  on $[\tau,T]\subseteq [0,T]$, such a solution is unique due to the local Lipschitz continuity of the vector field \cite{H69}. Taking transpose on both sides of \eqref{DE3_P} gives an ODE system for $\mbP_i^T$, $1\le i\le N$, which shows that $(\mbP_1^T, \cdots, \mbP_N^T)$ still satisfies \eqref{DE3_P}. So the ODE system \eqref{DE3_P}  guarantees  each $\mbP_i$ to be symmetric
\end{remark}

\begin{remark} \label{remark:PSr}
 If \eqref{DE3_P} has a unique
solution $(\mbP_1,\cdots ,\mbP_N)$ on $[0, T]$, then we can uniquely solve $(\mbS_1, \cdots, \mbS_N)$ and $(\mbr_1, \cdots, \mbr_N)$ by using  linear ODEs.
\end{remark}

For the $N$-player game, we consider closed-loop perfect state information, so that the state vector $X(t)$ is available to each player.

\begin{theorem}\label{theorem:Nash}
Suppose that \eqref{DE3_P} has a unique solution $(\mbP_1,\cdots ,\mbP_N)$ on $[0,T]$. Then we can uniquely solve \eqref{DE3_S}, \eqref{DE3_gamma}, and the game of $N$ players has a set of feedback Nash  strategies given by
$$
u_i=-R^{-1} \mbB_i^T (\mbP_i X(t) +\mbS_i), \quad 1\le i\le N.
$$
\end{theorem}
\proof This theorem follows the standard results in \cite[Theorem 6.16, Corollary 6.5]{BO99} .\endproof

By Theorem \ref{theorem:Nash},   the solution of the feedback Nash strategies with closed-loop perfect state information completely reduces to the study  of \eqref{DE3_P}.
For this reason, our subsequent analysis starts by analyzing \eqref{DE3_P}.

\section{Asymptotic Solvability}
\label{sec:asy}

For an $l\times m$ real matrix $Z =(z_{ij})_{i\le l, j\le m}$, denote the $l_1$-norm $\|Z\|_{l_1}= \sum_{i,j}|z_{ij}|$.
\begin{definition}\label{definition:as0}
The sequence of  Nash games \eqref{stateXi}-\eqref{costJi} has asymptotic solvability if there exists $N_0$ such that for all $N\ge N_0$,
 $(\mbP_1,\cdots,\mbP_N)$ in \eqref{DE3_P} has a solution on $[0,T]$ and
\begin{align}\label{main_conl1}
\sup_{N\ge N_0}\sup_{1\le i\le N, 0\leq t\leq T}  \|\mbP_i(t)\|_{l_1} <\infty.
\end{align}
\end{definition}

 Definition \ref{definition:as0} only involves the Riccati equations. This is sufficient due to Remark \ref{remark:PSr}. The boundedness condition
\eqref{main_conl1} is to impose certain regularity of the solutions, which is necessary for
studying the asymptotic behavior of the system when $N\to \infty$.

Let the $Nn\times Nn$ identity matrix be partitioned in the form:
\begin{align*}
I_{Nn}=\begin{bmatrix}
I_n   &0    &\cdots        &0    \\
0  &I_n    &\cdots         &0     \\
\vdots  &\vdots     &\ddots        &\vdots     \\
0  &0    &0         &I_n
\end{bmatrix}.
\end{align*}
For $1\le i\ne j\le N$,
exchanging the $i$th and $j$th rows of submatrices in $I_{Nn}$,  let $J_{ij}$ denote the resulting matrix. For instance, we have
\begin{align*}
J_{12}=\begin{bmatrix}
0   &I_n    &\cdots        &0    \\
I_n  &0    &\cdots         &0     \\
\vdots  &\vdots     &\ddots        &\vdots     \\
0  &0    &0         &I_n
\end{bmatrix}.
\end{align*}
It is easy to check that $J_{ij}^T=J_{ij}^{-1}=J_{ij}$.

\begin{theorem}\label{theorem:Prep3}
 We assume that \eqref{DE3_P} has  a solution
 $(\mbP_1(t), \cdots,
\mbP_N(t))$  on $[0,T]$. Then the following holds.

i) ${\mbP}_1(t)$  has the  representation
\begin{align}\label{P_matrix3}
{\mbP}_1(t)=\begin{bmatrix}
\Pi_1(t) &\Pi_2(t) &\Pi_2(t)&\cdots &\Pi_2(t) \\
\Pi_2^T(t) &\Pi_3(t) &\Pi_3(t)&\cdots &\Pi_3(t)\\
\Pi_2^T(t) &\Pi_3(t)&\Pi_3(t) &\cdots &\Pi_3(t)\\
\vdots          & \vdots        &  \vdots        &\ddots &\vdots \\
\Pi_2^T(t) &\Pi_3(t) &\Pi_3(t)&\cdots &\Pi_3(t)
\end{bmatrix},
\end{align}
where $\Pi_1$ and $\Pi_3$ are $n\times n$ symmetric matrices.

ii) For $i>1$,  $\mbP_i(t)= J_{1i}^T \mbP_1(t) J_{1i}$.
\end{theorem}

\proof See Appendix A. \endproof

By Theorem \ref{theorem:Prep3}, \eqref{main_conl1} is equivalent to the following condition:
\begin{align}\label{main_con2}
\sup_{N\ge N_0, 0\leq t\leq T} \left[|\Pi_1(t)|+N|\Pi_2(t)|+N^2|\Pi_3(t)|\right]<\infty.
\end{align}

We present some  continuous dependence result of parameterized  ODEs. This will play a  key role in establishing Theorem \ref{theorem:iff} below.
Consider
\begin{align}
&\dot{x}=f(t,x), \quad x(0)=z\in\mathbb{R}^K, \label{dot_x} \\
&\dot{y}=f(t,y)+g(\epsilon, t,y), \label{dot_y}
\end{align}
where $y(0)=z_{\epsilon}\in\mathbb{R}^K$, $ 0<\epsilon\leq 1.$

Let $\phi(t,x)=f(t,x),$ or $f(t,x)+g(\epsilon, t, x)$.

A1) $\sup_{\epsilon, 0\leq t\leq T} |f(t,0)|+|g(\epsilon, t, 0)|\leq C_1$.

A2) $\phi(\cdot, x)$ is Lebesgue measurable for each fixed $x\in\mathbb{R}^K$.

A3) For each $t\in [0, T]$, $\phi(t,x): \mathbb{R}^K\rightarrow\mathbb{R}^K$ is locally Lipschitz continuous in $x$, uniformly with respect to $(t,\epsilon)$, i.e., for any fixed $r>0$, and $x, y\in B_{r}(0)$ which is the open ball of radius $r$ centering $0$,
\begin{align*}
|\phi(t, x)-\phi(t, y)|\leq \mbox{Lip} (r)|x-y|,
\end{align*}
where $\mbox{Lip}(r)$ depends only on $r$, not on $\epsilon,t\in [0,T]$.

A4) For each fixed $r>0$,
\begin{align*}
\lim_{\epsilon\rightarrow 0}\sup_{0\leq t\leq T, y\in B_r(0)} |g(\epsilon, t, y)|=0, \quad \lim_{\epsilon\rightarrow 0}|z_{\epsilon}-z|=0.
\end{align*}

 If the solutions to \eqref{dot_x} and \eqref{dot_y}, denoted by $x^z(t)$ and $y^\epsilon(t)$,
 exist on $[0,T]$, they are unique by the local Lipschitz condition; for \eqref{dot_x} in this case denote
$\delta_\epsilon = \int_0^T |g(\epsilon, \tau, x^z(\tau))|d\tau$, which converges to $0$ as $\epsilon \to 0$ due to A4).
\begin{theorem} \label{theorem:depen}
i) If \eqref{dot_x} has a solution $x^z(\cdot)$ on $[0, T]$, then there exists $0<\bar{\epsilon}\leq 1$ such that for all $0<\epsilon\leq\bar{\epsilon}$, \eqref{dot_y} has a solution $y^\epsilon(\cdot)$ on $[0, T]$ and
\begin{align}
\sup_{0\le t\le T}|y^\epsilon (t)-x^z(t)|=O(|z_\epsilon-z|
+\delta_\epsilon). \label{oe}
\end{align}

ii) Suppose there exists a sequence $\{\epsilon_i, i\geq 1\}$ where $ 0<\epsilon_i\leq 1$ and $\lim_{i\to \infty}\epsilon_i= 0$ such that \eqref{dot_y} with $\epsilon=\epsilon_i$ has a solution $y^{\epsilon_i}$   on $[0, T]$   and
$\sup_{i\ge 1, 0\leq t\leq T} |y^{\epsilon_i}(t)|\leq C_2$
for some constant $C_2$. Then \eqref{dot_x} has a solution on $[0, T]$.
\end{theorem}

\proof  See Appendix B. \endproof

\begin{remark} \label{remark:tc}
If \eqref{dot_x} and \eqref{dot_y} are replaced by matrix ODEs and (or) a terminal condition at $T$ is used in each equation, the results in Theorem \ref{theorem:depen} still hold.
\end{remark}

 Let $$M=BR^{-1}B^T.$$ Before presenting further results, we introduce two Riccati ODEs:
\begin{align}\label{d11}
\begin{cases}
\dot{\Lambda}_1 =  \Lambda_1M \Lambda_1-(\Lambda_1A+A^T\Lambda_1)-Q, \\
\Lambda_1(T)=Q_f,
\end{cases}
\end{align}
and
\begin{align}\label{d21}
\begin{cases}
\dot{\Lambda_2} = \Lambda_1M \Lambda_2+ \Lambda_2 M \Lambda_1+ \Lambda_2M\Lambda_2 \\
                             \quad \qquad - (\Lambda_1G + \Lambda_2 (A+G) +A^T\Lambda_2) +Q\Gamma, \\
 \Lambda_2(T)=-Q_f\Gamma_f. 
 \end{cases}
\end{align}
Note that \eqref{d11} is the standard Riccati ODE in LQ optimal control and has a unique solution $\Lambda_1$ on $[0,T]$. Equation
\eqref{d21} is a non-symmetric Riccati ODE where $\Lambda_1$ is now treated as a known function. We state the main theorem on asymptotic solvability. The proof
is postponed to Section \ref{sec:proof}.

\begin{theorem}\label{theorem:iff}
The sequence of games  in  \eqref{stateXi}-\eqref{costJi} has asymptotic solvability if and only if   \eqref{d21} has a unique solution on $[0, T]$. \endproof
\end{theorem}

Our method of proving Theorem \ref{theorem:iff}  is to re-scale by
defining
\begin{align} \label{new_system_1}
&\Lambda_1^N=\Pi_1(t), \ \Lambda_2^N=N\Pi_2(t),
\ \Lambda_3^N=N^2\Pi_3(t),  
\end{align}
and examine their ODE system.
 We introduce the additional equation
\begin{align} \label{d3_1}
\begin{cases}
\dot{\Lambda}_3 = \Lambda_2^TM \Lambda_2 + \Lambda_3M \Lambda_1
+ \Lambda_1M \Lambda_3+\Lambda_3M \Lambda_2+\Lambda_2^TM\Lambda_3 \\
                             \qquad -
                             \left(\Lambda_2^T G+G^T\Lambda_2
                             +\Lambda_3( A+G)
                             +(A^T+G^T)\Lambda_3 \right)  \\
 \qquad - \Gamma^T Q\Gamma,  \\
  \Lambda_3(T)=\Gamma_f^T Q_f\Gamma_f.
\end{cases}
\end{align}
 Note that after \eqref{d11} and \eqref{d21} are solved on $[0,T]$ or otherwise on a maximal existence interval for the latter, \eqref{d3_1}  becomes  a linear ODE.

\begin{theorem} \label{theorem:PiLa}
Suppose \eqref{d21} has  a solution on $[0,T]$. Then we have
\begin{align}
\sup_{0\le t\le T}(|\Pi_1-\Lambda_1|  +| N\Pi_2 -
\Lambda_2|+|N^2\Pi_3-\Lambda_3|)=O(1/N).\nonumber
\end{align}
\end{theorem}

\proof The bound follows from Theorem \ref{theorem:depen} i) by use of
 $g_1, g_2, g_3$ and  the
initial conditions which appear in the equations of  $\Lambda_1^N$, $\Lambda_2^N$, $\Lambda_3^N$. \endproof

 \section{Proof of Theorem \ref{theorem:iff}}
 \label{sec:proof}



%

Note that $\Pi_3=\Pi_4$.
We rewrite the  system of \eqref{d1}, \eqref{d2} and
\eqref{d3}
 by use of a set of new variables
\begin{align*} 
&\Lambda_1^N=\Pi_1(t), \ \Lambda_2^N=N\Pi_2(t),
\ \Lambda_3^N=N^2\Pi_3(t).
\end{align*}
Here and hereafter $N$ is used as a superscript in various places.
This should be clear from the context.
We can determine functions $g_k$, $1\le k\le 3$, and obtain
\begin{align}
&\dot{\Lambda}_1^N = \Lambda_1^NM \Lambda_1^N-(\Lambda_1^NA+A^T\Lambda_1^N)
-Q \nonumber \\
&\qquad + g_1(1/N, \Lambda_1^N, \Lambda_2^N),\label{Nd1_1} \\
& \Lambda^N_1(T)=(I-\frac{\Gamma^T_f}{N})Q_f(I-\frac{\Gamma_f}{N}), \nonumber %
\end{align}
\begin{align}
&\dot{\Lambda}_2^N =  \Lambda_1^NM \Lambda_2^N+
\Lambda_2^N M \Lambda_1^N+
\Lambda_2^NM\Lambda_2^N \nonumber\\
&\qquad- (\Lambda_1^NG + \Lambda_2^N (G+A)
+A^T\Lambda_2^N) +Q\Gamma\nonumber \\
&\qquad + g_2(1/N,  \Lambda_2^N, \Lambda_3^N),  \label{Nd2_1} \\
&\Lambda^N_2(T)=-({I}-\frac{\Gamma_f^T}{N})Q_f\Gamma_f,  \nonumber
\end{align}
\begin{align}
&\dot{\Lambda}_3^N = (\Lambda_2^N)^TM \Lambda_2^N + \Lambda_3^NM \Lambda_1^N + \Lambda_1^NM \Lambda_3^N \nonumber \\
& \qquad +\Lambda_4^NM \Lambda_2^N+(\Lambda_2^N)^TM\Lambda_4^N\nonumber \\
                             &\qquad - \Big((\Lambda_2^N)^T G+G^T\Lambda_2^N+\Lambda_4^NG+ G^T\Lambda_4+\Lambda_3^N A+A^T\Lambda_3^N \Big) \nonumber \\
                             &\qquad  - \Gamma^T Q\Gamma \nonumber \\
                             &\qquad +  g_3(1/N, \Lambda_2^N, \Lambda_3^N),  \label{Nd3_1}  \\
& \Lambda_3^N(T)=\Gamma_f^T Q_f\Gamma_f . \nonumber
\end{align}

In particular, we can determine
\begin{align*}
g_1= &\frac{1}{N} (1-\frac{1}{N})( \Lambda_2^N M \Lambda_2^N
+(\Lambda_2^N)^T M (\Lambda_2^N)^T )\\
 & -\frac{1}{N} (\Lambda_1^N G + G^T\Lambda_1^N) \\
 &- \frac{1}{N} (1-\frac{1}{N}) (\Lambda_2^N G+ G^T (\Lambda_2^N)^T)\\
 &+ \frac{1}{N}(\Gamma^T Q+Q\Gamma)-\frac{1}{N^2}\Gamma^T Q \Gamma.
\end{align*}
The expressions of $g_2, g_3$ can be determined in a similar way and the detail is omitted here.


Note that if \eqref{d21} has a unique solution on $[0,T]$, we can uniquely solve $\Lambda_3$.
In view of $g_1, g_2, g_3$ and the terminal conditions in  \eqref{Nd1_1}-\eqref{Nd3_1},
by Theorem \ref{theorem:depen} and Remark \ref{remark:tc}, we obtain the desired result. \endproof

\section{Decentralized Control}
\label{sec:dec}

\begin{proposition}  Suppose that \eqref{DE3_P} has  a solution
 $(\mbP_1, \cdots,
\mbP_N)$  on $[0,T]$. Then $\mbS_i(t)$ in $\eqref{DE3_S}$ has the  form
\begin{align}\label{S_form}
\mbS_i(t)=[ \theta^T_2(t), \cdots,  \theta_1^T(t), \cdots,
\theta_2^T(t)]^T,
\end{align}
in which the $i$th sub-vector is $\theta_1(t)\in \mathbb{R}^n$ and the remaining sub-vectors are $\theta_2(t)\in \mathbb{R}^n$. Moreover, $\mbr_1=\cdots=\mbr_N$.
\end{proposition}

\proof
We can show that $$(J_{23}^T\mbS_1,\ J_{23}^T\mbS_3,\ J_{23}^T\mbS_2,\ J_{23}^T\mbS_4, \cdots, J_{23}^T\mbS_N)$$
satisfies \eqref{DE3_S} as $(\mbS_1, \cdots, \mbS_N)$ does.
Hence $\mbS_1= J_{23}^T \mbS_1$. We can further show $\mbS_1= J_{12}^T \mbS_2$. Following the argument in the proof of Lemma \ref{lem:Prep}, we obtain the representation \eqref{S_form}.

By \eqref{S_form}, we further obtain
\begin{align}
&\dot{\mbr}_i(t)= \theta_1^T M \theta_1 +2(N-1) \theta_2^TM \theta_1
 -\mbox{Tr}( D^T \Pi_1 D) \nonumber\\
 &\qquad\quad -(N-1) \mbox{Tr}( D^T \Pi_3 D) -\eta^T Q \eta \label{rie} \\
 & \mbr_i(T)= \eta_f^T Q_f \eta_f, \nonumber
\end{align}
for all $i$, so that $\mbr_1=\cdots =\mbr_N$,
\endproof

 Recalling $M=BR^{-1}B^T$, we derive
\begin{align}
&\dot{\theta_1}(t) = \Pi_1M\theta_1+(N-1)(\Pi_2M
\theta_1+\Pi_2^T M\theta_2)\nonumber\\
                            &\qquad\quad-\big(A^T+\frac{G^T}{N}\big)
                            \theta_1-
                            \frac{N-1}{N}G^T\theta_2 \nonumber \\
                           &\qquad\quad+\big(I-\frac{\Gamma^T}{N}
                           \big)Q\eta,  \label{theta_1} \\
                           & \theta_1(T)=-(I-\frac{\Gamma_f^T}{N}
                           )Q_f\eta_f,\nonumber \\
&\dot{\theta_2}(t) = \big(\Pi_2^T+(N-1)\Pi_3\big)M\theta_1\nonumber\\
                            &\qquad\quad+\big(\Pi_1+(N-2)\Pi_2^T\big)M
                            \theta_2-\frac{1}{N}G^T\theta_1 \nonumber \\
                            &\qquad\quad-\big(A^T+\frac{N-1}{N}G^T\big)\theta_2
                            -\frac{1}{N}\Gamma^TQ\eta, \label{theta_2}  \\
                            & \theta_2(T)=\frac{1}{N}\Gamma_f^T Q_f \eta_f. \nonumber
\end{align}

 Let
\begin{align*}
{\chi}_1^N(t)=\theta_1(t),\quad {\chi}_2^N(t)=N\theta_2(t).
\end{align*}
We may write
the ODEs of ${\chi}_1^N(t) $ and ${\chi}_2^N(t)$, which have the limiting form:
\begin{align}
&\dot{\chi}_1(t) =  (\Lambda_1M+\Lambda_2M-A^T){\chi_1} + Q\eta,\label{theta_1_new} \\
& {\chi_1}(T)=-Q_f\eta_f, \nonumber
\end{align}
and
\begin{align}
&\dot{\chi}_2(t) = ((\Lambda_2+\Lambda_4)M-G^T){\chi_1} \nonumber  \\
&\qquad\quad+((\Lambda_1+\Lambda_2)M
                        -(A^T+G^T)){\chi_2}-\Gamma^TQ\eta,
                        \label{theta_2_new}  \\
     &{\chi_2}(T)=\Gamma_f^T Q_f\eta_f. \nonumber
\end{align}

For \eqref{rie}, we have  the limiting form
\begin{align*}
&\dot\mbr(t) =\chi_1^T M\chi_1 + 2\chi_2^T M\chi_1 -\mbox{Tr}( D^T \Lambda_1 D)-\eta^T Q\eta,\\
& \mbr(T)=\eta_f^T Q_f \eta_f.
\end{align*}

\begin{proposition} \label{prop:thch}
If asymptotic solvability holds,
\begin{align}\label{sichi}
\sup_{0\le t\le T}(|\theta_1(t)-\chi_1(t)|+|N\theta_2(t)-\chi_2(t)|) =O(1/N).
\end{align}
\end{proposition}
\proof
Under asymptotic solvability, we uniquely solve $(\Lambda_1, \Lambda_2, \Lambda_3, \chi_1, \chi_2)$ on $[0,T]$. 
 We obtain \eqref{sichi} by writing the ODE system of
$(\Lambda_1^N, \Lambda_2^N, \Lambda_3^N, \chi_1^N, \chi_2^N)$ and next applying Theorem~\ref{theorem:depen}. \endproof

\subsection{Decentralized control and mean field dynamics}

By Theorem \ref{theorem:Nash}, the strategy of player ${\cal A}_i$ is
\begin{align}
u_i= -R^{-1} B^T\Big(\Pi_1(t)X_i+\Pi_2(t)\sum_{j\ne i} X_j +\theta_1(t)\Big).
\end{align}
The closed-loop equation of $X_i$ is  now given by
\begin{align}
&dX_i(t)= \Big(AX_i-M\Big(\Pi_1(t)X_i+\Pi_2(t)\sum_{j\ne i}
X_j \nonumber \\
& \qquad\qquad+\theta_1(t)\Big)
 +G X^{(N)}\Big)dt + D dW_i.
\end{align}
We introduce
$$
\frac{d\bar X}{dt}= \left(A-M (\Lambda_1 +\Lambda_2)+G\right)\bar X -M \chi_1(t),
$$
where $\Bar X(0)= x_0$, and further approximate
$X^{(N)}$ by $\bar X$. 
When $N\to \infty$, we obtain the decentralized control law
$$
u_i^d= -R^{-1} B^T \left(\Lambda _1(t) X_i+ \Lambda_2 (t) \bar X
+\chi_1(t)\right).  $$
 The next lemma provides an error estimate for the mean field approximation.
\begin{proposition}
Suppose  $E\sup_{i\ge 1}|X_i(0)|^2\le C$ for some fixed $C$ and
$\lim_{N\to \infty} \frac{1}{N} \sum_{i=1}^N EX_i(0)= x_0$.
Then
$$\sup_{0\le t\le T} E| X^{(N)}(t)- \bar X(t)|^2 =
O(|\frac{1}{N} \sum_{i=1}^N EX_i(0)-x_0|^2+1/N ). $$
\end{proposition}

\proof We first write the SDE for $X^{(N)}$ and find the explicit expression of $X^{(N)}(t)-\bar X(t)$. The proposition follows from elementary estimates by use of Theorem \ref{theorem:PiLa} and   Proposition \ref{prop:thch}.~\endproof


\section{An Example}
\label{sec:exa}

\begin{figure}[t]
\begin{center}
\begin{tabular}{c}
\psfig{file=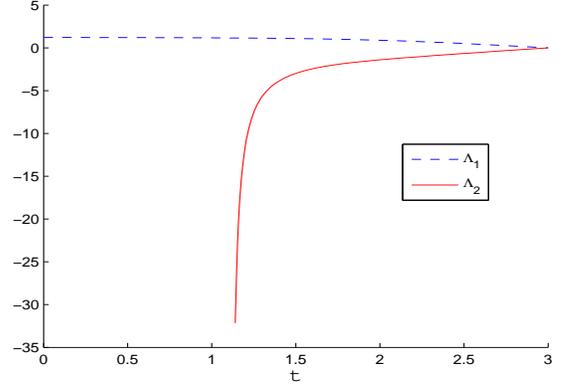, width=3.4in, height=2.2in}
\end{tabular}
\end{center}
\caption{ $\Lambda_2$ has a maximal existence interval small than $[0,T]$ }  \label{fig:bu}
\end{figure}

{\bf Example} \label{counterexample}
Take $A=0.2$, $B=G=Q=R=1$, $\Gamma=1.2$, $\Gamma_f=0$, $Q_f=0$, $T=3$. Consider the equation system \eqref{d11} and \eqref{d21}:
\begin{align*}
\dot{\Lambda}_1 &= \Lambda_1^2-0.4\Lambda_1-1, \quad \Lambda_1(T)=0, \\
\dot{\Lambda}_2 &= 2\Lambda_1\Lambda_2+\Lambda_2^2-(\Lambda_1+1.4\Lambda_2)+1.2, \quad \Lambda_2(T)=0.
\end{align*}
$\Lambda_1$ can be solved explicitly.
It is numerically illustrated in Fig. \ref{fig:bu} that $\Lambda_2$
does not have a solution on the whole interval $[0,T]$, implying no asymptotic solvability.

\section{Conclusion}
\label{sec:con}

This paper studies the asymptotic solvability problem for LQ
mean field games and obtains a necessary and sufficient condition via a
non-asymmetric Riccati ODE.
The re-scaling technique used in this paper can be extended to more general models in terms of dynamics and interaction patterns \cite{H10,HN11}. This will be reported in  our future work.

\section*{Appendix A: Proof of Theorem \ref{theorem:Prep3}}
\renewcommand{\theequation}{A.\arabic{equation}}
\setcounter{equation}{0}
\renewcommand{\thetheorem}{A.\arabic{theorem}}
\setcounter{theorem}{0}

We prove the following   lemma first.

\begin{lemma}
\label{lem:Prep}
 We assume that \eqref{DE3_P} has  a solution
 $(\mbP_1(t), \cdots,
\mbP_N(t))$  on $[0,T]$. Then the following holds.

i) ${\mbP}_1(t)$  has the  representation
\begin{align}\label{P_matrix}
{\mbP}_1(t)=\begin{bmatrix}
\Pi_1(t) &\Pi_2(t) &\Pi_2(t)&\cdots &\Pi_2(t) \\
\Pi_2^T(t) &\Pi_3(t) &\Pi_4(t)&\cdots &\Pi_4(t)\\
\Pi_2^T(t) &\Pi_4(t)&\Pi_3(t) &\cdots &\Pi_4(t)\\
\vdots          & \vdots        &  \vdots        &\ddots &\vdots \\
\Pi_2^T(t) &\Pi_4(t) &\Pi_4(t)&\cdots &\Pi_3(t)
\end{bmatrix}
\end{align}
where $\Pi_1$, $\Pi_3$ and $\Pi_4$ are $n\times n$ symmetric matrices.

ii) For $i>1$,  $\mbP_i(t)= J_{1i}^T \mbP_1(t) J_{1i}$.
\end{lemma}

\proof 
Step 1.  For $1\le i\le N$,  denote
$$\mbP_i= (\mbP_i^{jk})_{1\le j,k\le N}, $$
where each $\mbP_i^{jk}$ is an $n\times n$ matrix.
Define the new functions
$J_{23}^T \mbP_i J_{23}$, $i=1, \ldots, N$. By elementary calculations, we  see that
$$(J_{23}^T \mbP_1 J_{23}, \ J_{23}^T \mbP_3 J_{23},\  J_{23}^T \mbP_2 J_{23},\ J_{23}^T \mbP_4 J_{23}, \cdots, J_{23}^T \mbP_N J_{23})$$
satisfies \eqref{DE3_P} together with its terminal condition as $(\mbP_1(t), \cdots,
\mbP_N(t))$ does. Hence
$
\mbP_1= J_{23}^T \mbP_1 J_{23},
$
which implies
\begin{align}
\mbP_1^{12}=\mbP_1^{13},\quad  \mbP_1^{22}=P_1^{33},
\quad \mbP_1^{23}=\mbP_1^{32}.\label{P123}
\end{align}
Repeating the above by using $J_{2k}$, $k\ge 4$, in place of $J_{23}$, we obtain
$$
\mbP^{12}_1=\mbP_1^{13}=\cdots=\mbP_1^{1N}, \qquad \mbP_1^{22}=\mbP_1^{33}=\cdots =\mbP_1^{NN}.
$$
We similarly obtain $\mbP_1= J_{24}^T \mbP_1 J_{24}$, and this gives
$$
\mbP^{23}_1=\mbP^{24}_1.
$$
Repeating the similar argument,
we can check all other remaining off-diagonal submatrices. Since $\mbP_1$ is symmetric (also see Remark \ref{remark:P}), $(\mbP_1^{23})^T= \mbP_1^{32}$, $\mbP_1^{23}$  is symmetric by \eqref{P123}. By the above method we can show that all off-diagonal submatrices on neither the first row nor the first column are identical and symmetric. Therefore we obtain the representation of $\mbP_1$.

Step 2.  We can verify that both
$$(J_{12}^T \mbP_2 J_{12}, \ J_{12}^T \mbP_1 J_{12},\  J_{12}^T \mbP_3 J_{12}, \cdots, J_{12}^T \mbP_N J_{12})$$ and
$(\mbP_1(t), \cdots,
\mbP_N(t))$ satisfy \eqref{DE3_P}. Hence $\mbP_2=J_{12}^T \mbP_1 J_{12}$.
All other cases can be similarly checked.
\endproof

{\it Proof of Theorem \ref{theorem:Prep3}:}
By Lemma \ref{lem:Prep}, we have
\begin{align}
&\dot{\Pi_1}(t)  =   \Pi_1 M\Pi_1 
                          +(N-1)(\Pi_2M\Pi_2
                          +\Pi_2^T M\Pi_2^T) \nonumber \\
                          &\quad\qquad - \Big(\Pi_1 (A+\frac{G}{N})
                          +(A^T+\frac{G^T}{N})\Pi_1 \Big)\nonumber \\
                          &\quad\qquad-(1-\frac{1}{N}) (\Pi_2G
                          +G^T\Pi_2^T) \nonumber \\
                          &\quad\qquad -({I}-\frac{\Gamma^T}{N})Q({I}-\frac{\Gamma}{N}),
                        \label{d1}   \\
                  &\Pi_1(T)=({I}-\frac{\Gamma_f^T}{N})
                 Q_f({I}-\frac{\Gamma_f}{N}),\nonumber
              \end{align}
and
\begin{align}
&\dot{\Pi_2}(t) =  \Pi_1 M\Pi_2+\Pi_2M\Pi_1
+\Pi_2^TM\Pi_3\nonumber \\
             &\quad\qquad+(N-2)(\Pi_2M\Pi_2
           +\Pi_2^TM\Pi_4) \nonumber \\
             &\quad\qquad-\Big(\Pi_1 \frac{G}{N}+\frac{G^T}{N}\Pi_3
          +\frac{N-2}{N}G^T\Pi_4 \nonumber \\
       &\quad\qquad+\Pi_2(A+\frac{N-1}{N}G)
      +(A^T+\frac{G^T}{N})\Pi_2\Big) \nonumber \\
& \quad\qquad+ ({I}-\frac{\Gamma^T}{N})Q\frac{\Gamma}{N},\label{d2} \\
&\Pi_2(T)=-({I}-\frac{\Gamma_f^T}{N})Q_f\frac{\Gamma_f}{N},\nonumber
\end{align}
and
\begin{align}
&\dot{\Pi_3}(t) = \Pi_2^TM\Pi_2+\Pi_3M\Pi_1
+\Pi_1 M\Pi_3\nonumber \\
                       &\quad\qquad+(N-2)(\Pi_4M\Pi_2
                       +\Pi_2^TM\Pi_4) \nonumber\\
                       &\quad\qquad - \Big( \frac{1}{N}(\Pi_2^TG
                       +G^T\Pi_2)
                        \nonumber \\
                       &\qquad\qquad+\Pi_3(A+\frac{G}{N})
                       +(A^T+\frac{G^T}{N})\Pi_3\nonumber \\
                       & \qquad\qquad +\frac{N-2}{N}(\Pi_4G +
                       G^T\Pi_4)
                        \Big)\nonumber \\
&\quad\qquad - \frac{\Gamma^T}{N}Q\frac{\Gamma}{N}, \label{d3} \\
& \Pi_3(T)=\frac{\Gamma_f^T}{N}Q_f\frac{\Gamma_f}{N},\nonumber
\end{align}
and
\begin{align}
&\dot{\Pi_4}(t) =\Pi_2^TM\Pi_2
+\Pi_4 M\Pi_1
+\Pi_1 M\Pi_4 \nonumber \\
                       &\quad\qquad+\Pi_3M\Pi_2
                       +\Pi_2^TM\Pi_3 \nonumber \\
                       &\quad\qquad +(N-3)(\Pi_4 M\Pi_2
                       +\Pi_2^TM\Pi_4 ) \nonumber \\
                       &\quad\qquad - \Big( \frac{1}{N}(\Pi_2^T G
                       +G^T\Pi_2
                       +\Pi_3G
                       +G^T\Pi_3) \nonumber \\
                       &\quad\qquad+ \Pi_4(A+\frac{N-2}{N}G)
                       +(A^T+\frac{N-2}{N}G^T)\Pi_4 \Big)\nonumber \\
&\quad\qquad - \frac{\Gamma^T}{N}Q\frac{\Gamma}{N}, \label{d4} \\
& \Pi_4(T)=\frac{\Gamma_f^T}{N}Q_f\frac{\Gamma_f}{N}.\nonumber
\end{align}

Then we can further show that $\Pi_3-\Pi_4$ satisfies a linear ODE when $\Pi_1$ and $\Pi_2$ are fixed and that $\Pi_3(T)-\Pi_4(T)=0$. This gives $\Pi_3=\Pi_4 $ on $[0, T]$. \endproof

\section*{Appendix B: Proof of Theorem \ref{theorem:depen}}
\renewcommand{\theequation}{B.\arabic{equation}}
\setcounter{equation}{0}
\renewcommand{\thetheorem}{B.\arabic{theorem}}
\setcounter{theorem}{0}

\begin{proof}
i) Let $x^z(t)$ be the solution of \eqref{dot_x} on $[0,T]$, and we can find a constant $C_z$ such that $\sup_{0\leq t\leq T} |x^z(t)|\leq C_z$, and $\sup_{0<\epsilon\le 1}|z_\epsilon|\le C_z$.
Fix the open ball $B_{2C_z}(0)$. For $x, y\in B_{2C_z}(0)$ and $t\in [0,T]$, we have
\begin{align*}
|\phi(t, x)-\phi(t, y)|\leq \mbox{Lip} (2C_z) |x-y|.
\end{align*}

For each $\epsilon\le 1$,
by A1)-A3), \eqref{dot_y} has a solution $y^{\epsilon}(t)$  defined either (a) for all $t\in[0, T]$ or (b) on a maximal interval $[0, t_{\max})$ for some $0<t_{\max}<T$.

Below we  show that for all small $\epsilon$, (b) does not occur.
We prove by contradiction. Suppose
for any small $\epsilon_0>0$, there exists $0<\epsilon<\epsilon_0$ such that (b) occurs with the corresponding 
 $0< t_{\max} <T$.
Since $[0, t_{\max})$ is the maximal existence interval,
we have $\lim_{t\uparrow t_{\max}}|y^\epsilon(t)|=\infty$ \cite{H69}.
Therefore for some $0<t_m<t_{\max}$,
\begin{align}
y^{\epsilon}(t_m)\in \partial B_{2C_z}(0) \label{con_3},
\end{align}
and
\begin{align}
y^{\epsilon}(t)\in B_{2C_z}(0), \quad \forall 0\le t<t_m. \label{con_4}
\end{align}
For $t< t_{\max}$,
 we have
\begin{align*}
&y^{\epsilon}(t)-x^z(t)=z_{\epsilon}-z\\
&+\int_0^t
\Big[f(\tau, y^{\epsilon}(\tau))+g(\epsilon,
\tau, y^{\epsilon}(\tau))-f(\tau, x^z(\tau))\Big]d\tau.
\end{align*}
Denote $
\zeta(\tau) =f(\tau, y^{\epsilon}(\tau))+g(\epsilon,
\tau, y^{\epsilon}(\tau))-f(\tau, x^z(\tau))$ and
it follows that
\begin{align*}
|\zeta(\tau)| &= |\zeta(\tau)-g(\epsilon,
\tau, x^z(\tau))+ g(\epsilon,
\tau, x^z(\tau))|\\
& \le \mbox{Lip}(2C_z) |y^\epsilon(\tau) - x^z(\tau)|
+|g(\epsilon,\tau, x^z(\tau))|.
\end{align*}
Now for $0\le t<t_m$,
\begin{align*}
|y^{\epsilon}(t)-x^z(t)| \leq&\ |z_{\epsilon}-z|+\delta_{\epsilon}\\
 &+\int_0^t \mbox{Lip} (2C_z)|y^{\epsilon}(\tau)-x^z(\tau)|d\tau.
\end{align*}
Note that $\delta_{\epsilon}=\int_0^T|g\big(\epsilon, \tau, x^z(\tau)\big)|d\tau\rightarrow 0$ as $\epsilon\rightarrow 0$.
By Gronwall's lemma,
\begin{align*}
|y^{\epsilon}(t)-x^z(t)|\leq (\delta_{\epsilon}+|z_{\epsilon}-z|)e^{\mbox{Lip}(2C_z)t}
\end{align*}
for all $t\leq t_m$.
We can find $\bar{\epsilon}>0$ such that for all $\epsilon\leq \bar{\epsilon}$,
\begin{align*}
(\delta_{\epsilon}+|z_{\epsilon}-z|)e^{\mbox{Lip}(2C_z)T}<\frac{C_z}{3}.
\end{align*}
Then for all $0\leq t\leq t_m$,
$y^{\epsilon}(t)\in B_{3C_z/2}(0)$,
which is a contradiction to \eqref{con_3}.
We conclude for all $0<\epsilon \leq \bar{\epsilon}$, $y^{\epsilon}$ is defined on $[0, T]$. Next, \eqref{oe} follows readily.

ii) We have
\begin{align}
y^{\epsilon_i}(t)=z_{\epsilon_i}+\int_0^t \Big[f\big(\tau, y^{\epsilon_i}(\tau)\big)+g\big(\epsilon, \tau, y^{\epsilon_i}(\tau)\big)\Big]d\tau, \label{con_5}
\end{align}
and
\begin{align}
&|f\big(\tau, y^{\epsilon_i}(\tau)\big)
+g\big(\epsilon, \tau, y^{\epsilon_i}(\tau)\big)|\nonumber  \\
&\leq \mbox{Lip} (C_2)|y^{\epsilon_i}(\tau)|
 +|f(\tau, 0)+g(\epsilon, \tau, 0)| \nonumber \\
&\leq \mbox{Lip} (C_2)|y^{\epsilon_i}(\tau)|+C_1 \nonumber \\
& \leq \mbox{Lip}(C_2) C_2 +C_1 \label{con_6},
\end{align}
where $C_1$ is given in A1).

By \eqref{con_5}-\eqref{con_6}, the functions $\{y^{\epsilon_i}(\cdot), i\geq 1\}$ are uniformly bounded and equicontinuous. By Ascoli's lemma, there exists a subsequence $\{y^{\epsilon_{i_j}}(\cdot), j=1, 2, 3,\cdots\}$ such that $y^{\epsilon_{i_j}}$ converges to $y^*\in C\big([0, T], \mathbb{R}^K\big)$ uniformly on $[0, T]$.
Hence,
\begin{align*}
y^*(t)=z+\int_0^t f\big(\tau, y^*(\tau)\big)+g\big(\epsilon, \tau, y^*(\tau)\big)d\tau
\end{align*}
for all $t\in [0, T]$. So $\eqref{dot_x}$ has a solution.
\end{proof}

The  proof in part i) follows the method in
 \cite[sec. 2.4]{P96} and \cite[pp. 486]{S98}.

\end{document}